\documentclass{amsart}

\usepackage{amsmath,amssymb,amsthm,latexsym}

\newtheorem{theorem}{Theorem}
\newcommand{\bt}{\begin{theorem}}
\newcommand{\et}{\end{theorem}}
\newtheorem{lemma}{Lemma}
\newcommand{\bl}{\begin{lemma}}
\newcommand{\el}{\end{lemma}}
\newtheorem{corollary}{Corollary}
\newcommand{\bc}{\begin{corollary}}
\newcommand{\ec}{\end{corollary}}
\newcommand{\bconj}{\begin{conjecture}}
\newcommand{\econj}{\end{conjecture}}
\newtheorem{problem}{Problem}
\newcommand{\bprob}{\begin{problem}}
\newcommand{\eprob}{\end{problem}}
\newcommand{\beq}{\begin{equation}}
\newcommand{\eeq}{\end{equation}}
\newcommand{\benum}{\begin{enumerate}}
\newcommand{\eenum}{\end{enumerate}}
\newcommand{\N}{\ensuremath{ \mathbf N }}
\newcommand{\Z}{\ensuremath{\mathbf Z}}

\newcommand{\mcf}{\ensuremath{ \mathcal F}}

\newcommand{\mcr}{\ensuremath{ \mathcal R}}

\newcommand{\bmat}{\left(\begin{matrix}}
\newcommand{\emat}{\end{matrix}\right)}
\newcommand{\bsmallmat}{\left(\begin{smallmatrix}}
\newcommand{\esmallmat}{\end{smallmatrix}\right)}

\DeclareMathOperator{\qqand}{\qquad\text{and}\qquad}

\title[Sizes of sumsets]{Problems in additive number theory, VI: \\
Sizes of sumsets of finite sets}
\author{Melvyn B.  Nathanson}
\address{Department of Mathematics\\Lehman College (CUNY)\\Bronx, NY 10468}
\email{melvyn.nathanson@lehman.cuny.edu}

\date{\today}

\thanks{Supported in part by  PSC-CUNY Research Award Program grant 66197-00 54.}

\begin{document}

\begin{abstract}
In the study of sums of finite sets of integers, most attention has been paid to sets with small sumsets 
(Freiman's theorem and related work) and to sets with large sumsets (Sidon sets and $B_h$-sets).  
This paper focuses on the full range of sizes of $h$-fold sums of a set of $k$ integers.   
Many new results and open problems are presented.  
\end{abstract}

\maketitle

\section{A remark of Erd\H os and Szemer\' edi} 

Let $\N = \{1,2,3,\ldots\}$ be the set of positive integers,  $\N_0 = \{0,1,2,3,\ldots\}$ 
 the set of nonnegative integers, and  \Z\ the set of all integers. 
In this paper, the letters $h,i,j,k,\ell,m$, and $n$ always denote positive integers.  

For real numbers $u$ and $v$, the \emph{interval of integers}  from $u$ to $v$ is the set 
\[
[u,v] = \{c \in \Z: u \leq c \leq v \}.
\]
The half-open, half-closed intervals $(u,v]$ and $[u,v)$ are defined similarly.
 
Let $G$ be an additive abelian group  
and let $\mcf^*(G)$ be the set of all nonempty finite subsets of $G$. 
For every integer $h$ and every nonempty subset $A$ of $G$, 
the $h$-fold \emph{sumset}  of $A$, denoted $hA$,
 is the set of all sums of $h$ not necessarily distinct elements of $A$:
\[
hA = \underbrace{A+\cdots + A}_{\text{$h$ summands}}
= \left\{ a_1+\cdots + a_h: a_i \in A \text{ for all } i \in [1,h] \right\}. 
\]
The $h$-fold \emph{restricted sumset}  of $A$, denoted $\widehat{hA}$, 
is the set of all sums of $h$  distinct elements of $A$:
\begin{align*}
\widehat{hA} 
& = \left\{ a_1+\cdots + a_h: a_i \in A \text{ for all } i \in [1,h] 
\text{ and } a_i \neq a_j \text{ for all } i \neq j \right\} \\ 
& = \left\{ \sum_{s\in S}s :   S \subseteq A \text{ and } |S| = h \right\}
\end{align*}
where $|S|$ denotes the size (or cardinality) of the set $S$. 
The sumset $hA$ and the restricted sumset $\widehat{hA}$ are finite if $A$ is finite.  

In a paper on the sum-product conjecture, 
 Erd\H os and  Szemer\' edi~\cite{erdo-szem83} wrote:
\begin{quotation}
Let $2n-1 \leq t \leq \frac{n^2+n}{2}$.  It is easy to see that one can find a sequence 
of integers $a_1 < \ldots < a_n$ so that there should be exactly $t$ distinct integers 
in the sequence $a_i+a_j, 1 \leq i \leq j \leq n$.  
\end{quotation}
This is a statement about the sizes of 2-fold sumsets of finite sets of integers.  
The purpose of this paper is to investigate  the sizes of $h$-fold sumsets 
and $h$-fold restricted sumsets of finite subsets of the integers and, 
more generally, of additive abelian groups for all $h \geq 2$. 

The size and structure of  small $h$-fold sumsets and  small $h$-fold 
restricted sumsets  for finite subsets 
of  the integers and of finite fields have been extensively investigated 
(cf.~\cite{abbo24}--\cite{nath96bb}), 
but there has been considerably less study of the set of \emph{all} 
$h$-fold sumsets and $h$-fold restricted sumsets 
of finite subsets of an abelian group $G$.  
We consider the \emph{range of cardinalities of $h$-fold sumsets} 
and \emph{$h$-fold restricted sumsets} 
of sets of size $k$ in $G$.  
For all positive integers $h$ and $k$,   define   
\[
\mcr_G(h,k) = \left\{ \left| hA \right|: A \subseteq G \text{ and } |A| = k \right\}.   
 \]
Thus, 
\[
\widehat{\mcr}_G(h,k) = \left\{ \left|\widehat{hA} \right|: A \subseteq G \text{ and } |A| = k \right\}.
 \]

There is a simple symmetry for restricted sumsets.  
Let $|A| = k$ and $\Lambda = \sum_{a\in A}a$.  Let $h \in [1,k-1]$.
If $S \subseteq A$ and $|S| = h$, then  
\[
b = \sum_{a\in S}a \in \widehat{hA}. 
\]
We have  $|A\setminus S| = k-h \in [1,h-1]$ and 
\[
\sum_{a\in A\setminus S} a = \Lambda - b \in \widehat{(k-h)A}. 
\] 
Thus, $\Lambda - \widehat{hA} = \widehat{(k-h)A}$ and 
$\left| \widehat{hA} \right| = \left| \widehat{(k-h)A} \right|$.  
It follows that   
\beq                 \label{sizes:symmetry}
\widehat{\mcr}_G(h,k) = \widehat{\mcr}_G(k-h,k).
\eeq

We also observe that $h>k$ implies $\widehat{hA} = \emptyset$ 
 and $\widehat{\mcr}_G(h,k) = \{0\}$.

\bprob
An open problem in additive number theory is to compute the range of sumset sizes  $\mcr_{G}(h,k)$ 
and $\widehat{\mcr}_{G}(h,k)$.
\eprob


\section{Upper and lower bounds for sumset sizes in groups}

There are simple best possible   upper and lower bounds for the sets 
$\mcr_G(h,k)$ and $\widehat{\mcr}_G(h,k)$ in an additive abelian group $G$.

\bt                                                   \label{sizes:theorem:simple} 
Let $G$ be an additive abelian group.  For all  integers $h\geq 2$ and $k\geq 2$, 
there are best possible upper bounds 
\beq                                                      \label{sizes:simple-upper} 
\max \mcr_G(h,k) \leq \binom{k+h-1}{h}
\eeq 
and 
\beq                                                      \label{sizes:simple-upper-hat} 
\max \widehat{\mcr}_G(h,k) \leq \binom{k}{h}  
\eeq
and best possible lower bounds 
\beq                                                      \label{sizes:simple-lower} 
\min \mcr_G(h,k) \geq k 
\eeq 
and 
\beq                                                     \label{sizes:simple-lower-hat} 
\min \widehat{\mcr}_G(h,k) \geq \max(k-h+1,h+1).  
\eeq 
Moreover,  for all $k \geq 1$, 
\beq                                                    \label{sizes:simple-lower-1} 
\mcr_G(1,k) = \widehat{\mcr}_G(1,k) = \widehat{\mcr}_G(k-1,k)  = \{k\}
\eeq
and,  for all $h \geq 1$, 
\beq                                                    \label{sizes:simple-lower-hat-1} 
\mcr_G(h,1) = \{1\}   
\qqand  
 \widehat{\mcr}_G(h,1) = \begin{cases}
 1 & \text{if $h=1$}\\
0 & \text{if $h \geq 2$.}
 \end{cases} 
\eeq 
\et

\begin{proof} 
Let  $A$ be a subset of $G$  with $|A| = k$. 
For all positive integers $h$, we have the combinatorial  upper bounds
\[
|hA| \leq \binom{k+h-1}{h}
\qqand 
\left| \widehat{hA} \right| \leq  \binom{k}{h}. 
\]
These upper bounds are best possible.  
For example, in the additive group \Z, 
for every integer $g \geq h+1$, every set $A$ of $k$ distinct powers of $g$ 
gives the upper bounds in inequalities~\eqref{sizes:simple-upper}    
and~\eqref{sizes:simple-upper-hat} 
(cf. Nathanson~\cite{nath96bb}, Theorems 1.3 and 1.9). 
Another example: For all $i \in [1,k]$, let $m_i > h$ and let $G = \bigoplus_{i=1}^k \Z/m_i\Z$.
 The set 
 \[
 A = \{(0,\ldots,0,1 +m_i\Z, 0,\ldots,0) :i \in [1,k]\}
 \]
 also gives the upper bounds in~\eqref{sizes:simple-upper}  and~\eqref{sizes:simple-upper-hat} . 

For the lower bounds for sumsets $hA$, we observe that if $G$ is a group 
and if $A \subseteq G$ with $|A|=k$, 
then, for all $a, a_0 \in A$, we have  
\[
 (h-1)a_0 + a \in hA  
\] 
and so   
\[               
k = \left| \{(h-1)a_0 + a: a \in A\} \right| \leq |hA|.
\]
This proves~\eqref{sizes:simple-lower}.
For every  finite subset  $A$ of the group $G$ with $|A| = k$, 
we have $|hA| = k$ for some $h \geq 2$ if and only if $A$ 
is a coset of a subgroup of $G$ of order $k$. 
(This is a known result in group theory.  For completeness, 
 a proof is in Appendix~\ref{sizes:appendix:hA=A}.) 
It follows that the lower bound~\eqref{sizes:simple-lower} is best possible. 

Next, we compute   lower bounds for restricted sumsets $\widehat{hA}$ for $h \leq k$.  

If $h \in [1,k-1]$ and $a_1,\ldots, a_{h-1}$ are distinct elements of $A$, then 
\[
\left\{ a_1+\cdots + a_{h-1} + a: a \in A \setminus \{a_1,\ldots, a_{h-1}\} \right\} 
\subseteq \widehat{hA} 
\]
and so   
\[
\left| \widehat{hA} \right| \geq k-h+1 
\]
and 
\[
\min \widehat{\mcr}_G(h,k) \geq k-h+1. 
\]
The symmetry~\eqref{sizes:symmetry} implies 
\[
\min \widehat{\mcr}_G(h,k) = \min \widehat{\mcr}_G(k-h,k) \geq k-(k-h)+1 = h+1  
\]
and so 
\[
\min \widehat{\mcr}_G(h,k) \geq \max(k-h+1,h+1). 
\]
This proves~\eqref{sizes:simple-lower-hat}.

Identities~\eqref{sizes:simple-lower-1} and~\eqref{sizes:simple-lower-hat-1} 
follow immediately from the definitions of the sumset size sets.  
This completes the proof. 
\end{proof}

Consider   the group \Z.  
It follows from Theorem~\ref{sizes:theorem:simple} that, for all positive integers 
$h$ and $k$,  the sumset size set  $\mcr_{\Z}(h,k)$ 
is finite, and so, for all $t \in \mcr_{\Z}(h,k)$, 
there is a set $A$ with $\min(A) = 0$ such that $|A| = k$ and $|hA|=t$.  
Thus, there is an integer $N = N(h,k)$ such that, for all $t \in \mcr_{\Z}(h,k)$, 
there is a set $A \subseteq [0,N]$ such that $|A| = k$ and $|hA|=t$.  

\bprob
Compute an integer $N(h,k)$ with this property.  Compute the least integer with this property. 
\eprob
\bprob
Compute the analogous integers $\widehat{N}(h,k)$ for the restricted sumset size set $\widehat{\mcr}_{\Z}(h,k)$.
\eprob

\bprob
Let $G$ be an ordered abelian group or a group of unbounded exponent.  
Let $h$, $k$, and $t$ be positive integers 
such that 
\[
hk-h+1 \leq t \leq \binom{k+h-1}{h}.  
\]
Does there exist a subset $A$ of $G$ with $|A| = k$ and $|hA| = t$? 
This is a decision problem in the complexity class NP. 
Is this  problem NP-complete? 
\eprob

\bprob
Let $G$ be an ordered abelian group or a group of unbounded exponent.  
Let $h$, $k$, and $t$ be positive integers 
such that $3 \leq h \leq k$ and 
\[
hk-h^2+1 \leq t \leq \binom{k}{h}.  
\]
Does there exist a subset $A$ of $G$ with $|A| = k$ and $\left| \widehat{hA} \right| = t$? 
This is a decision problem in the complexity class NP. 
Is this  problem NP-complete? 
\eprob


\section{Ordered groups and groups of unbounded exponent}      \label{sizes:section:integers} 

Our principal interest is  sumset sizes of sets of integers, 
but many results about sumset sizes of sets of integers depend only on the order
structure of the integers and extend to all ordered abelian groups.  
Indeed, many proofs extend to groups that simply contain elements 
of arbitrarily large finite order or elements of infinite order .   

A multiplicative group $G$, not necessarily abelian, is   \emph{ordered} 
with respect to a total order $<$ on $G$ 
if, for all $x,y,z \in G$, the inequality $x < y$ implies $xz < yz$ and $zx < zy$.  
A group $G$ is \emph{orderable} if $G$ is an ordered group with 
respect to some total order on $G$.
If $x \in G$ and $x \neq e$, then either $e < x$ or $e > x$. 
If $e < x$, then $e < x < x^2 < \cdots$.   
If $x < e$, then $e > x > x^2 > \cdots$ .  In both cases, $x$ has infinite order.  
Thus, every ordered group is torsion-free, but not 
every torsion-free group is orderable. 
Torsion-free abelian groups and  torsion-free nilpotent groups are orderable (Glass~\cite{glas99}).  

An additive abelian group $G$ is \emph{ordered} with respect to a total order $<$ on $G$ 
if, for all $x,y,z \in G$, the inequality $x < y$ implies $x+z < y+z$.

An additive abelian group $G$ has \emph{bounded exponent} if 
there is a positive integer $m$ such that $mx=0$ for all $x \in G$.  
An abelian group $G$ has \emph{exponent $m$} if $G$ has bounded exponent and 
$m$ is the smallest positive integer such that $mx=0$ for all $x \in G$.  
If $G$ is a torsion group such that every element of $G$ has order at most $M$, 
then the exponent of $G$ is bounded and is at most the least common multiple 
of the orders of the elements of $G$.  
The group $G$ has \emph{unbounded exponent} if it does not have bounded exponent. 
A group that is not a torsion group contains an element of infinite order, 
and so has unbounded exponent.  Every ordered group is torsion-free 
and so is of unbounded exponent.  The group $G = \bigoplus_{n=1}^{\infty} \Z/n\Z$ 
is a torsion group of unbounded exponent. 

In the integers and in all ordered abelian groups and  abelian groups of unbounded exponent, 
we have lower bounds for the sizes of $h$-fold sumsets 
and restricted sumsets that are better than the inequalities in Theorem~\ref{sizes:theorem:simple}. 

\bt         \label{sizes:theorem:LowerBound-G} 
Let $G$ be an ordered abelian group. 
For all positive integers $h$ and $k$,   
there are the lower bounds 
\beq                              \label{sizes:hA-ordered-size-lower}
\min \mcr_G(h,k)  =  hk-h+1 
\eeq
and 
\beq                              \label{sizes:hA-ordered-hat-size-lower}
\min \widehat{\mcr}_G(h,k)  = hk-h^2+1.
\eeq
and the upper bounds 
\beq                                                      \label{sizes:simple-upper-1} 
\max \mcr_G(h,k) = \binom{k+h-1}{h}
\eeq 
and 
\beq                                                      \label{sizes:simple-upper-hat-1} 
\max \widehat{\mcr}_G(h,k) = \binom{k}{h}  
\eeq

\et

\begin{proof}
Let $G$ be an ordered group and let $A = \{a_i:i=1,\ldots, k\}$ 
be a $k$-element subset of $G$ with 
\[
a_1< \cdots < a_i < a_{i+1} < \cdots < a_k.
\]
For all $i \in [1, k-1]$ we have the strictly increasing sequences 
\[
ha_1 < \cdots < ha_i < ha_{i+1} < \cdots < ha_k 
\]
and  
\[
ha_i  < (h-1) a_i +  a_{i+1} <  (h-2) a_i + 2 a_{i+1} < \cdots < a_i  + (h-1) a_{i+1}<  h a_{i+1}.
\]
These inequalities give the following lower bound for $\left| hA \right|$: 
\[
 \left| hA \right| \geq h(k-1)+1 = hk-h+1.
\]

To obtain a lower bound for $\left| \widehat{hA} \right|$, we consider, 
for all $i \in [ 1,k-h]$  and $j \in [0,h]$, 
the strictly increasing sequence of $h+1$ elements 
$s_{i,j}$ of the group that are obtained from the $(h+1)$-fold sum
\[
a_i  + a_{i+1} + \cdots + a_{i+h-j-1}  + a_{i+h-j}  + a_{i+h-j+1} + \cdots + a_{i+h}
\]
by deleting the summand $a_{i+h-j}$.  Thus,  
\[
s_{i,j} = a_i  + a_{i+1} + \cdots + a_{i+h-j-1}  + a_{i+h-j+1} + \cdots + a_{i+h}
\in \widehat{hA}. 
\] 
In particular, 
\[
s_{i,0} =   a_i  + a_{i+1} +  \cdots + a_{i+h-1}.  
\]
The inequalities  
\[
s_{1,0} <  \cdots < s_{i,0} <  s_{i+1,0}  < \cdots < s_{k-h,0}  < s_{k-h+1,0} = \sum_{j=1}^{h} a_{k-h+j}
\]
and 
\[
s_{i,0} < \cdots < s_{i,j-1} <   s_{i,j} < \cdots < s_{i,h} = s_{i+1,0}   
\]
imply the lower bound 
\[
 \left| \widehat{hA} \right| \geq h(k-h)+1 =  hk-h^2+1.
\]
 
Every arithmetic progression of length $k$ in an ordered group $G$ 
gives the lower bounds in inequalities~\eqref{sizes:hA-ordered-size-lower} 
and~\eqref{sizes:hA-ordered-hat-size-lower} and so these lower bounds are best possible.  
This proves the lower bounds.  
The upper bounds follow from Theorem~\ref{sizes:theorem:simple} and its proof. 
\end{proof}

The following simple but important remark allows us to transport results 
about sumsets of the integers 
to sumsets of any additive abelian group of unbounded exponent.  

\bt              \label{sizes:theorem:unboundedTransfer}
If the abelian group $G$ has unbounded exponent, then 
\[
\mcr_{\Z}(h,k)  \subseteq \mcr_G(h,k) 
\]
and 
\[
\widehat{ \mcr}_{\Z}(h,k)  \subseteq \widehat{\mcr}_G(h,k).
\]
\et

\begin{proof}
Let $A$ be a finite set of integers with $|A| = k$ and let $h \geq 2$.    
By translation, we can assume that $\min A = 0$.  
If $a = \max A $, then  $ha = \max hA $.  Because $G$ has unbounded exponent, 
there exists $x \in G$ of order greater than $ha$,  
and so the set $A_G = \{ax:a \in A\}$ is a subset of $G$ of size $k$. 
Then 
\[
hA_G= \{ bx:b\in hA\} \qqand \widehat{hA_G} = \{cx:c \in \widehat{hA} \}.
\]
Because $x$ has order greater than $ha$, we have $|hA_G| = |hA|$ 
and $\left|\widehat{hA_G}\right| = \left|\widehat{hA}\right|$. 
Thus, $t \in \mcr_{\Z}(h,k)$ implies $t \in \mcr_G(h,k)$ and 
$t \in\widehat{ \mcr}_{\Z}(h,k)$ implies $t \in\widehat{\mcr}_G(h,k)$. 
This completes the proof.  
\end{proof}

It is not necessarily true that $\mcr_{\Z}(h,k) = \mcr_G(h,k)$ 
and $\widehat{\mcr}_{\Z}(h,k) = \widehat{\mcr}_G(h,k)$.  
For example, let $h \geq 2$ and $k \geq 2$.  
Let
\[
m = hk-h 
\]
and 
\[
G = (\Z/m\Z )\oplus \Z. 
\]
The group $G$ has unbounded exponent and contains the subset 
\[
A = \left\{ (i,0): i \in[0,k-1] \right\}
\] 
of size $k$, but 
\begin{align*}
hA & = \{(j,0):j \in [0, hk-h] \} \\ 
&  =  \{(j,0):j \in [0, m-1] \} \\ 
& = \Z/m\Z \oplus \{0\}. 
\end{align*}
Thus, $|hA| = m = hk-h$ and $m \in \mcr_G(h,k)$, but  
$\min \mcr_{\Z}(h,k) = hk-h+1$ by Theorem~\ref{sizes:theorem:LowerBound-G}  
and so $m \notin \mcr_{\Z}(h,k)$.  
 
 Similarly, let 
 \[
 m = hk-h^2
 \]
 and 
 \[
 G = (\Z/m\Z )\oplus \Z. 
 \]
The group $G$ has unbounded exponent and contains the subset  
 \[
A = \left\{ (i,0): i \in[0,k-1] \right\}
\] 
of size $k$.
If $ A_0 = [0,k-1]$, then the sumset $h A_0$ is the interval of integers 
\[
\widehat{h A_0} =  \left[ \frac{h(h-1)}{2}, hk- \frac{h(h+1)}{2} \right]  
\]
of cardinality  
\[
\left| \widehat{h A_0} \right| =  hk - \frac{h(h+1)}{2}  - \frac{h(h-1)}{2} + 1 =  hk-h^2 +1 = m+1.
\] 
It follows that 
\begin{align*}
\widehat{hA} 
= \left\{ (j,0):j \in h A_0 \right\} 
=  \Z/m\Z \oplus \{0\} 
\end{align*}
and so $m = |\widehat{hA}|  \in \widehat{\mcr}_G(h,k)$. 
By Theorem~\ref{sizes:theorem:LowerBound-G},   
\[
m = hk-h^2 <  hk-h^2+1 = \min \widehat{\mcr}_{\Z}(h,k)
\] 
and so $m \notin \widehat{\mcr}_{\Z}(h,k)$.  
Thus, $\widehat{\mcr}_{\Z}(h,k)$ is a proper subset of $\widehat{\mcr}_{G}(h,k)$. 


\bt                     \label{sizes:theorem:LowerBound-G-unbounded}
Let $G$ be an abelian group of unbounded exponent.  
For all positive integers $h$ and $k$,   
\beq                              \label{sizes:hA-ordered-size-lower-exp}
\min \mcr_G(h,k) \leq  hk-h+1 
\eeq
and 
\beq                              \label{sizes:hA-ordered-hat-size-lower-exp}
\min \widehat{\mcr}_G(h,k)  \leq   hk-h^2+1.
\eeq
\et

\begin{proof}
This follows immediately from Theorems~\ref{sizes:theorem:LowerBound-G}
and~\ref{sizes:theorem:unboundedTransfer}, applied to the ordered group \Z.  
This completes the proof.  
\end{proof}

\section{Sizes of 2-fold sumsets and 2-fold restricted sumsets}     \label{sizes:section:h=2}

Let $G$ be an additive abelian group  of unbounded exponent.  
We shall prove that, for every positive integer $k$ and every integer $t$ 
such that 
\[
2k-1 \leq t \leq \frac{k^2+k}{2}
\]
there is a $k$-element subset $A$ of $G$ with sumset size $|2A| = t$.  
We also prove that, for every positive integer $k$ and every integer $t$ such that 
 \[
 2k-3 \leq t \leq \frac{k^2-k}{2}
 \] 
 there is a $k$-element subset $A$ of $G$ with restricted sumset size
  $\left|\widehat{2A}\right| = t$.  
By Theorem~\ref{sizes:theorem:unboundedTransfer}, it suffices to prove these results for 
the integers. 

In every ordered abelian group, the range of sizes of 2-fold sumsets 
and the range of sizes of 2-fold restricted sumsets are intervals of consecutive integers.   
There is an exponential time algorithm to compute a set $A$ with $|A| = k$ 
and $|2A| = t$ for all $t \in \mcr_{\Z}(2,k)$ and 
to compute a set $A$ with $|A| = k$ 
and $|\widehat{2A}| = t$ all $t \in  \widehat{\mcr}_{\Z}(h,k)$. 
  
These results do not necessarily hold in groups of bounded exponent.  
Little is known about sumset sizes in these groups, even for the group 
$G = \bigoplus_{i=1}^{\infty} \Z/2\Z$ of exponent 2.

\bt                        \label{sizes:theorem:RZ2k}
For every ordered abelian group $G$ and for all positive integers $k$, 
\[
\mcr_{\Z}(2,k) = \mcr_{G}(2,k) = \left[ 2k-1, \frac{k^2 +k}{2} \right].
\]
Moreover, for all $t \in \mcr_{\Z}(2,k)$, there exists a set $A \subseteq \left[0,2^k -1 \right]$ 
such that $|A| = k$ and $\left|\widehat{2A} \right| = t$.
\et 

\begin{proof}
By Theorems~\ref{sizes:theorem:LowerBound-G} and~\ref{sizes:theorem:unboundedTransfer},
we have  
\[
\mcr_{\Z}(2,k)  \subseteq \mcr_G(2,k) 
\]
and 
\[
\min \mcr_{\Z}(2,k)  = \min \mcr_G(2,k) = 2k-1 
\]
and 
\[
\max \mcr_{\Z}(2,k)  = \max \mcr_G(2,k) = \frac{k^2 +k}{2} .  
\]
Proving that the set  $\mcr_{\Z}(2,k)$ is an interval of consecutive integers 
suffices to prove the Theorem.  

If $A$ is a set  of $k$ integers 
with $|2A| = t$, then $k$ and $t$  satisfy inequality
\[
2k-1 \leq t \leq \frac{k^2 +k}{2}. 
\]
For $k = 1$, choose $A = \{0\}$.  We have $2A = \{0\}$ and $|2A| = 1$, 
and so $\mcr_{G}(2,1) = [1,1]$. 

For $k = 2$, choose  $A = \{0,1\}$.  We have $2A = \{0,1,2\}$  and $|2A| = 3$, 
and so $\mcr_{G}(2,2) = [3,3]$.

For $k = 3$, choose $A _1= \{0,1,2\}$ and  $A_2 = \{0,1,3\}$.   
We have $2A_1 = \{0,1,2,3,4 \}$  with $|2A_2| = 5$ 
and  $2A_2 = \{0,1,2,3,4,6\}$ with $|2A| = 6$, 
and so $\mcr_{G}(2,3) = [5,6]$.

Thus, Theorem~\ref{sizes:theorem:RZ2k} is true for $k=1$, $2$, and $3$.  

For all $k \geq 4$, let $b$ be an integer such that 
\[
k-1 \leq b \leq 2k-3 
\]
and let 
\[
A = [0, k-2] \cup \{b\}.
\]
Then $|A| = k$ and 
\[
A \subseteq [0, b] \subseteq [0,2k-3]  \subseteq \left[0,2^k - 1 \right].
\] 
Because $k-2 < b \leq 2k-3$, we have  
\begin{align*}
2A 
& = [0, 2k-4] \cup [b, b+ k -2\} \cup \{ 2b\} \\
& = [0,b+k-2 ] \cup \{2b\} 
\end{align*}
and so 
\[
|2A| = b+k.
\]
This proves the theorem for 
\beq        \label{sizes:Z-1}
2k-1 \leq t \leq 3k-3. 
\eeq

We use induction on $k$ to prove the theorem for 
\[
3k-3 \leq t \leq \frac{k^2+k}{2}.
\]
Let $k \geq 4$ and assume that Theorem~\ref{sizes:theorem:RZ2k} is true for $k-1$. 
Let $s$ be an integer such that 
\[
2k-3 = 2(k-1)-1 \leq s \leq \frac{(k-1)^2 +k-1}{2} = \frac{k^2 - k}{2}.
\] 
By the induction hypothesis, there exists a set $B$ of   integers 
such that 
\[
|B|= k-1 \qqand B\subseteq \left[0,2^{k-1}-1 \right] 
\]
and 
\[
|2B| = s. 
\]
Let 
\[
a_k = 2\max(B)+1 = \max(2B)+1. 
\]
We have $a_k \notin B$ and 
\[
a_k \leq 2\left(2^{k-1}-1\right)+1 = 2^k -1.
\] 
Let $A = B \cup \{a_k\}$.
Then 
\[
|A| = k \qqand A \subseteq \left[0,a_k\right] \subseteq \left[0,2^k -1 \right] 
\]
and 
\[
2A = 2B \cup \left(a_k+B \right) \cup \{2a_k\}.
\]
Because 
\[
\max(2B) = a_k-1 < a_k \leq  \min(a_k+B)
\]
the sets $2B$, $a_k+B$, and $\{2a_k\}$ are pairwise disjoint.  
It follows that  
\[
|2A| = |2B| + |B| + 1= s + k 
\]
and 
\[
3k-3 \leq s + k \leq \frac{k^2 + k}{2}.
\]
Thus,  for all $t$ such that 
\beq        \label{sizes:Z-2}
3k-3 \leq t \leq \frac{k^2 +k}{2} 
\eeq
there exists a set  $A  \subseteq \left[0, 2^k -1 \right]$ with $|A| = k$ and $|2A| = t$.
Combining inequalities~\eqref{sizes:Z-1} and~\eqref{sizes:Z-2} completes the proof.  
\end{proof}

There is an analogous result for restricted sumsets.

\bt                        \label{sizes:theorem:RZ2k-hat}
For every ordered abelian group $G$ and for all   integers $k \geq 2$, 
\[
\widehat{\mcr}_{\Z}(2,k) = \widehat{\mcr}_{G}(2,k) = \left[ 2k-3, \frac{k^2 - k}{2} \right].
\]
Moreover, for all $t \in\widehat{\mcr}_{\Z}(2,k)$, there exists a set $A \subseteq \left[0,2^{k-2} \right]$ 
such that $|A| = k$ and $\left|\widehat{2A} \right| = t$.
\et 

\begin{proof}
By Theorems~\ref{sizes:theorem:LowerBound-G} and~\ref{sizes:theorem:unboundedTransfer},
we have  
\[
\widehat{\mcr}_{\Z}(2,k)  \subseteq \widehat{\mcr}_G(2,k) 
\]
and 
\[
\min \widehat{\mcr}_{\Z}(2,k)  = \min \widehat{\mcr}_G(2,k) = 2k-3  
\]
and 
\[
\max \widehat{\mcr}_{\Z}(2,k)  = \max \widehat{\mcr}_G(2,k) = \frac{k^2 - k}{2} .  
\]
Proving that the set  $\widehat{\mcr}_{\Z}(2,k)$ is an interval of consecutive integers 
suffices to prove the Theorem.  

If $A$ is a set  of $k$ integers 
with $|\widehat{2A}| = t$, then $k$ and $t$  satisfy inequality
\[
2k-3 \leq t \leq \frac{k^2 -k}{2}. 
\]
For $k = 2$, choose $A = \{0,1\}$.  We have $\widehat{2A} = \{1\}$  and $|\widehat{2A}| = 1$, 
and so $\widehat{\mcr}_{G}(2,2) = [1,1]$. 

For $k = 3$, choose  $A = \{0,1,2\}$.  We have $\widehat{2A} = \{1,2,3 \}$  and $|\widehat{2A}| = 3$, 
and so $\widehat{\mcr}_{G}(2,3) = [3,3]$. 

For $k = 4$, choose  $A_1 = \{0,1,2,3\}$ and $A_2 = \{0,1,2,4\}$.  
We have $\widehat{2A_1} = \{1,2,3,4,5\}$  with  $|\widehat{2A_1}| = 5$, 
and  $\widehat{2A_2} = \{1,2,3,4,5,6\}$  with  $|\widehat{2A}| = 6$, 
and so $\widehat{\mcr}_{G}(2,4) = [5,6]$. 

Thus, Theorem~\ref{sizes:theorem:RZ2k-hat} is true for $k=2$, $3$, and $4$.  

For all $k \geq 5$, let $b$ be an integer such that 
\[
k-1 \leq b \leq 2k-4 
\]
and let 
\[
A = [0, k-2] \cup \{ b\}.
\]
Then $|A| = k$ and 
\[
A \subseteq [0, b] \subseteq [0,2k-4]  \subseteq \left[0, 2^{k-2} \right].
\] 
Because $b \leq 2k-4$, we have  
\begin{align*}
\widehat{2A}  
& = [1, 2k-5] \cup [b, b+k-2]  \\
& = [1, b+k-2] 
\end{align*}
and so 
\[
\left| \widehat{2A} \right|= b+k-2.
\]
This proves the theorem for 
\beq        \label{sizes:Z-1-hat}
2k-3 \leq t \leq 3k-6. 
\eeq 

We use induction on $k$ to prove the theorem for 
\[
3k-5 \leq t \leq \frac{k^2-k}{2}.
\]
Let $k \geq 5$ and assume that Theorem~\ref{sizes:theorem:RZ2k-hat} is true for $k-1$. 
Let $s$ be an integer such that 
\[
2k-5 = 2(k-1)-3 \leq s \leq \frac{(k-1)^2 - (k-1)}{2} = \frac{(k-1)(k-2)}{2}.
\] 
By the induction hypothesis, there exists a set $B$ of   integers 
such that 
\[
|B|= k-1 \qqand B\subseteq \left[0,2^{k-3} \right] 
\]
and 
\[
\left| \widehat{2B} \right| = s. 
\]
Let 
\[
a_k = \max \left( \widehat{2B} \right) +1 = \max(B) + \left(\max(B)-1\right) + 1 = 2\max(B). 
\]
We have $\max(B) \leq 2^{k-3}$ and so 
\[
a_k \leq 2^{k-2}.
\] 
Let $A = B \cup \{a_k\}$.
Then 
\[
|A| = k \qqand A \subseteq \left[0,a_k\right] \subseteq \left[0, 2^{k-2} \right] 
\]
and 
\[
\widehat{2A} = \widehat{2B} \cup \left(a_k+B \right).
\]
Because 
\[
\max \left( \widehat{2B} \right) = a_k-1 < a_k \leq \min(a_k+B)
\]
the sets $ \widehat{2B}$ and $a_k+B$  are pairwise disjoint and so 
\[
\left| \widehat{2A} \right| = \left| \widehat{2B} \right|  + |B|  = s + k - 1
\]
and 
\begin{align*}
3k-6 
& =(2k-5)+k -1\leq s + k -1 \\
& \leq  \frac{(k-1)(k-2)}{2} + (k-1) \\ 
& = \frac{k^2 - k}{2}.
\end{align*}
Thus,  for all $t$ such that 
\beq        \label{sizes:Z-2-hat}
3k-6 \leq t \leq \frac{k^2 - k}{2} 
\eeq 
there exists a set  $A  \subseteq \left[0, 2^{k-2} \right]$ with $|A| = k$ and $|2A| = t$.
Combining inequalities~\eqref{sizes:Z-1-hat} and~\eqref{sizes:Z-2-hat} completes the proof.  
\end{proof}


\bprob 
By Theorem~\ref{sizes:theorem:RZ2k}, for all $t \in \mcr_{\Z}(2,k) =  \left[ 2k-1, (k^2 +k)/2 \right]$, 
there is an exponential time algorithm to compute 
a set $A$ of integers with $|A| = k$ and $|2A| = t$: 
Examine all $k$-element subsets of the interval 
$[0,2^k -1]$.    Is there a polynomial time algorithm?  
Is there a polynomial bound on the set $A$?  
Equivalently, do there exist  integers $n$ and $c $ such that, 
for all $t \in \mcr_{\Z}(2,k)$, there is a set $A \subseteq \left[0, ck^n\right]$ 
with $|A| = k$ and $\left|\ 2A\right| = t$?
\eprob

\bprob
Do there exist  integers $n$ and $c $ such that, such that, 
 for all $t \in \widehat{\mcr} (2,k) $, there is a set $A \subseteq \left[0, ck^n\right]$ 
with $|A| = k$ and $\left|\widehat{2A} \right| = t$?
\eprob


\section{Missing sumset sizes}       \label{sizes:section:missing}
Let $G$ be an ordered additive abelian group.  
By Theorem~\ref{sizes:theorem:RZ2k}, the sumset size set $\mcr_{G}(2,k)$ 
is an interval  of consecutive integers.  
We also have $\mcr_{G}(1,k) = \{k\}$, $\mcr_G(h,1) = \{ 1\}$ 
and $\mcr_G(h,2) = \{h+1\}$ for all positive integers $h$ and $k$.    

It is natural to ask: Is every sumset size set  $\mcr_{G}(h,k)$ an interval? 
We prove that the answer is `no' for all $h \geq 3$ and $k \geq 3$.  
 The smallest example is the case $h=k=3$.  We have 
\[
\mcr_G(3,3) \subseteq \{7,8,9,10\}.
\]
For $A_1=\{0,1,2\}$, we have  $3A_1=[0,6]$ and $|3A_1| = 7$. \\
For $A_2=\{0,1,3\}$, we have  $3A_2=[0,7] \cup \{9\}$ and $|3A_2| = 9$.\\
For $A_3=\{0,1,4\}$, we have  $3A_3=[0,6] \cup \{8,9,12\}$ and $|3A_3| = 10$.\\
Thus,
\[
 \{7, 9,10\} \subseteq  \mcr_G(3,3).
\]
It follows from Theorem~\ref{sizes:theorem:Missing} that size 8 is missing from the set $\mcr_G(3,3)$. 
Equivalently, there exists no set of three integers whose 3-fold sum has size 8, and 
\[
\mcr_G(3,3) = \{7, 9,10\}.
\]  
Why?

\bt             \label{sizes:theorem:Missing}
Let $G$ be an ordered additive abelian group.  
For all $h \geq 3$ and $k \geq 3$, 
\[
 hk -h + 1 = \min\left( \mcr_{G}(h,k) \right) 
\]
but
\[
hk -h + 2 \notin \mcr_{G}(h,k). 
\]
Thus, the sumset size  set $\mcr_{G}(h,k)$ 
is not an interval of consecutive integers.  
\et

\begin{proof}
By Theorem~\ref{sizes:theorem:LowerBound-G}, we have 
\begin{align*}
\left\{hk-h+1, \binom{k+h-1}{h}\right\} & \subseteq \mcr_{G}(h,k) \\
&   \subseteq  \left[ hk-h+1, \binom{k+h-1}{h}  \right]. 
\end{align*} 
Let $A = \{a_1,\ldots, a_k\}$ be a $k$-element subset of 
the ordered group $G$, where 
\[
a_1 < a_2 < \cdots < a_k.
\]
For all $i \in [1,k-1]$ and $j \in [0,h-1]$, we have 
\begin{align*}
ha_i   < (h & -1)a_i    + a_{i+1}   < \cdots < (h-j)a_i + ja_{i+1} \\
& < \cdots 
<  a_i +  (h-1)a_{i+1} 
<  ha_{i+1}   
\end{align*} 
and so the sumset $hA$ contains the subset
\[
B = \{ha_k\} \cup \bigcup_{i=1}^{k-1} \left\{ (h-j)a_i + ja_{i+1} : j \in [0,h-1] \right\}. 
\]
It follows that 
\[
|hA| \geq |B| =h(k-1)+1 = hk-h+1. 
\]
We shall prove that if $hA \neq B$, then $|hA| \geq |B|+2 = hk-k+3$.  
Equivalently, we shall prove that if $hA \neq B$, then the sumset $hA$ 
must contain at least 
two elements that do not belong to $B$. 

For $i \in \{ 2,\ldots, k-1\}$ we have the inequalities 
\[
 a_{i-1} + (h-1)a_i  <  ha_i  < (h-1)a_i  + a_{i+1}  
\]
and 
\[
 a_{i-1} + (h-1)a_i  <   a_{i-1} + (h-2)a_i + a_{i+1} < (h-1)a_i  + a_{i+1}
\]
If   $hA = B$, then, for all $i \in [2,k-1]$, we have   
\[
 a_{i-1} + (h-2)a_i + a_{i+1} = ha_i
\]
and so 
\[
a_i -  a_{i-1} = a_{i+1}-a_i. 
\]  
Thus, $hA = B$ only if $A$ is a $k$-term arithmetic progression in $G$.  
(Conversely, if $A$ is a $k$-term arithmetic progression, then $hA = B$.)   

If $hA \neq B$, then  $A$ is not a $k$-term arithmetic progression in $G$ 
and there exists $j \in  [2,k-1]$ such that 
\beq             \label{sizes:badTriple} 
2a_j \neq  a_{j-1} +   a_{j+1}.  
\eeq
The following are three consecutive elements of the set $B$:   
\[ 
a_{j-1} + (h-1) a_j  < ha_j   <  (h-1)a_j + a_{j+1}. 
\]
We  have $a_{j-1} + (h-2)a_j + a_{j+1} \in hA$ and the inequality 
\[
a_{j-1} + (h-1) a_j  < a_{j-1} + (h-2)a_j + a_{j+1}   <  (h-1)a_j + a_{j+1}. 
\]
Inequality~\eqref{sizes:badTriple} implies 
\[
ha_j \neq  a_{j-1} + (h-2)a_j + a_{j+1} 
\] 
and so 
\[
a_{j-1} + (h-2)a_j + a_{j+1}  \in (hA) \setminus B. 
\]

We consider two cases.  In the first case,  
\beq                     \label{sizes:case 1}
a_{j-1} + (h-2)a_j + a_{j+1} < ha_j 
\eeq 
and, in the second case,  
\beq                     \label{sizes:case 2}
ha_j < a_{j-1} + (h-2)a_j + a_{j+1} .
\eeq

We  have $2a_{j-1} + (h-3)a_j + a_{j+1} \in hA$.  
In the first case,  inequality~\eqref{sizes:case 1} implies 
\begin{align*}
2a_{j-1} + (h-2) a_j  & <  2a_{j-1} + (h-3)a_j + a_{j+1}  \\
&  <  a_{j-1} + (h-2)a_j + a_{j+1} \\
&  < ha_j.  
\end{align*}
Also, 
\begin{align*}
2a_{j-1} + (h-2) a_j & <   a_{j-1} + (h-1)a_j \\
&   <  a_{j-1} + (h-2)a_j + a_{j+1} \\
&  < ha_j.
\end{align*}
It follows from~\eqref{sizes:badTriple} that 
\[
a_{j-1} + (h-1)a_j \neq  2a_{j-1} + (h-3)a_j + a_{j+1} 
\] 
and so 
\[
2a_{j-1} + (h-3)a_j + a_{j+1} \in (hA) \setminus B. 
\]
Thus, in the first case, $|hA\setminus  B| \geq 2$. 

We  have $a_{j-1} + (h-3)a_j + 2a_{j+1} \in hA$.  
In the second case,  inequality~\eqref{sizes:case 2}  implies  
\begin{align*}
ha_j  & <  a_{j-1} + (h-2)a_j + a_{j+1}  \\
&  <  a_{j-1} + (h-3)a_j + 2a_{j+1} \\
&  <  (h-2)a_j + 2a_{j+1}.  
\end{align*}
Also, 
\begin{align*}
ha_j & <  a_{j-1} + (h-2)a_j + a_{j+1} \\
& <  (h-1)a_j    + a_{j+1}  \\ 
& < (h-2)a_j + 2a_{j+1}.  
\end{align*}
It follows from~\eqref{sizes:badTriple} that 
\[
 (h-1)a_j    + a_{j+1}  \neq  a_{j-1} + (h-3)a_j + 2a_{j+1} 
\] 
and so 
\[
a_{j-1} + (h-3)a_j + 2a_{j+1}  \in (hA) \setminus B. 
\]
Thus, in the second case, $|hA\setminus  B| \geq 2$. 
This completes the proof. 
\end{proof}

Applying Theorem~\ref{sizes:theorem:Missing} with $ h=3$, we obtain 
\[
 \mcr_{G}(3,k) \subseteq  \{3k-2\} \cup \left[3k, \frac{k(k+1)(k+2)}{6} \right].
\]
We have proved that $\mcr_G(3,3)    = \{7\} \cup [9,10]$.

\bt
\begin{align}
\mcr_G(3,4) &  = \{10\} \cup [12,20]       \label{sizes:lemma:h=3,k=4}    \\
\mcr_G(3,5)  & = \{ 13\} \cup [15,35].               \label{sizes:lemma:h=3,k=5}  
\end{align}
\et

 \begin{proof}
 The proof is by the explicit construction of sets of integers.  
 
Let $P_{1,b} = \{0,1,2,b\}$ and  $P_{2,b} = \{0,1,3,b\}$.  Then 
\[
\left\{  \left| 3P_{1,b} \right|  : b \in [3,7]  \right\} = \{10,12,14,15,16 \} 
\]
\[
\left\{  \left| 3P_{2,b}   \right| : b \in [4,10]  \right\} = \{13,14,15,17,18,19 \} 
\]
and 
\[
 \left| 3\{ 1,4,16, 64 \}  \right|  =  \left| 3\{ 0,1,4,13 \}  \right|  = 20 .
\]
This proves~\eqref{sizes:lemma:h=3,k=4}.

Let $Q_{1,b} = \{0,1,2,3,b\}$, $Q_{2,b} = \{0,1,3,4,b\}$,   $Q_{3,b} = \{0,1,4,5,b\}$, 
$Q_{4,b} = \{0,1,5,7,b\}$, and $Q_{5,b} = \{ 0,1,5,8,19+3b \}$. 
Then   
\begin{align*} 
\left\{  \left|  3Q_{1,b}  \right|  : b \in [4,10]  \right\} & = \{ 13,15,17,19,20,21,22 \} \\
\left\{  \left| 3Q_{2,b} \right|  : b \in [5,13]  \right\} & = \{16,18,19,21,23,24,25,26,27 \} \\
\left\{  \left| 3Q_{3,b} \right|  : b \in [11,16]  \right\} & = \{26,27,28,29,30 \} \\
\left\{  \left| 3Q_{4,b} \right|  : b \in [19,22]  \right\} & = \{31,32,33 \} \\
\left\{  \left| 3Q_{5,b}  \}  \right|  : b \in [0,2]  \right\} & = \{32,33,34 \}
\end{align*} 
and 
\[
 \left| 3\{ 1,4,16, 64,256 \}  \right|  = 35.
\]
This proves~\eqref{sizes:lemma:h=3,k=5}.
\end{proof}

\bprob 
Compute $\mcr_{G}(3,k)$ for all $k \geq 6$. 
\eprob

\bprob
A general problem is to compute $\mcr_{\Z}(h,k)$ for fixed $h$ 
as $k$ increases.
\eprob


\section{The sumset size set $\mcr_{\Z}(h,3)$}
Let $G$ be an ordered group and let $A = \{a,b\}$ be a subset of $G$ of size 2. 
If $a < b$, then $hA = \{(h-i)a+ib : i \in [0,h]\}$  and  
\[
ha < (h-1)a+b <   \cdots < (h-i)a+ib < \cdots < a + (h-1)b < hb. 
\] 
It follows that  
\[
\mcr_G(h,2) = \{h+1\}.
\]
We would like to understand the sumset size set $\mcr_G(h,3)$. 
Theorem~\ref{sizes:theorem:Missing} with $ k=3$ gives  
\[
\mcr_{G}(h,3) \subseteq 
\{2h+1\} \cup \left[ 2h+3,  \binom{h+2}{2} \right]. 
\]

Let $G = \Z$.  Computing sumset sizes of random sets of integers of size 3, 
we obtain the left sides of the following relations: 
\begin{alignat*}{4}
                                       \{3\}  &   =       && \mcr_{\Z}(1,3) &  =                & & [3,3]       \\  
                                       \{5,6\}  &   =       && \mcr_{\Z}(2,3) &  =                & & [5,6]           \\  
                 \{7, 9,10\}   &   =             && \mcr_{\Z}(3,3) & \subseteq          & &  [7,10]             \\  
                    \{ 9,12,14,15\} &  \subseteq &&\mcr_{\Z}(4,3)& \subseteq & & \left[ 9, 15\right] \\  
          \{ 11, 15, 18, 20, 21\} & \subseteq && \mcr_{\Z}(5,3) &\subseteq & &  \left[ 11,  21 \right] \\
   \{ 13, 18, 22, 25, 27, 28 \}  &\subseteq && \mcr_{\Z}(6,3)& \subseteq & & \left[ 13,  28 \right] \\
      \{ 15, 21,26,30,33,35, 36 \}  &\subseteq && \mcr_{\Z}(7,3)& \subseteq & & \left[ 15,  36 \right] \\
         \{ 17, 24,30,35,39,42,44,45 \}  &\subseteq && \mcr_{\Z}(8,3)& \subseteq & &\left[ 17,  45 \right] \\
 \{ 19, 27,34,40,45, 49,52,54,55 \}  &\subseteq && \mcr_{\Z}(9,3)& \subseteq & &\left[ 19,  55 \right] \\
\{ 21, 30,38,45,51,56,60,63,65,66\}  &\subseteq && \mcr_{\Z}(10,3)& \subseteq & &  \left[ 21,66 \right]\\
\{ 23,33,42,50,57,63,68,72,75,77,78\}  &\subseteq && \mcr_{\Z}(11,3)& \subseteq & & \left[ 23,78 \right]\\
\{ 25, 36,46,55,63,70,76, 81, 85,88,90,91\}  &\subseteq && \mcr_{\Z}(12,3)& \subseteq & &  \left[ 25,91 \right] \\
\{ 27,39,50,60,69,77,84,90,95,99,102,104,105\}  &\subseteq && \mcr_{\Z}(13,3)& \subseteq & & \left[ 27,105 \right]
\end{alignat*}
Searching for a pattern in these examples, we find that the left side is 
\[
\left\{  \binom{h+2}{2} - \binom{t}{2} :  t \in [1,h] \right\}. 
\]
This suggests the following result.  

\bt
For all positive integers $h$, 
\beq                  \label{sizes:h-3} 
\mcr_{\Z}(h,3) = \left\{  \binom{h+2}{2} - \binom{t}{2} :  t \in [1,h] \right\}. 
\eeq  
\et

\begin{proof}
Every set of three integers is affinely equivalent to a set of the form 
$A = \{0,a,b\}$ with $0 < a < b$ and $\gcd(a,b) = 1$.  
Thus, it suffices to consider only sets of this form.  

Let $x_0,y_0,x_1,y_1$ be  integers with 
\[
x_0,x_1\in [0,b-1].
\]
If 
\[
x_0a+y_0 b =  x_1a+y_1b 
\]
and $y_0 > y_1$, then 
\[
 0 < b \leq  (y_0-y_1)b = (x_1-x_0)a 
\]
and $x_1 > x_0$.  
Because $a$ divides $ (y_0-y_1)b$ and $\gcd(a,b)=1$, we have 
$y_0-y_1=ra$ for some positive integer $r$ and so 
\[
rab = (x_1-x_0)a > 0.  
\]
Then 
\[
rb = x_1-x_0  
\]
and  
\[
b \leq rb =  x_1-x_0   \leq b-1
\]
which is absurd.  Therefore, $y_0=y_1$ and $x_0 = x_1$, and the representation 
of an integer in the form $xa+yb$ with $x\in [0,b-1]$ is unique.   

Let $n \in hA$.  There exist nonnegative integers $x$ and $y$ such that $x+y \leq h$
and 
\[
n = xa+yb.
\]
By the division algorithm, there are integers $q$ and $x_0$ 
such that $x = qb+x_0$ with $q \geq 0$ and $x_0 \in [0,b-1]$. 
Let $y_0 = y+qa$.  Then 
\[
n = (x-qb)a + (y+qa)b = x_0 a+ y_0 b
\]
where  
\begin{align*}
x_0+y_0 & = (x-qb)+(y+qa) = (x+y)-q(b-a) \\
&  \leq x+y-q  \leq x+y \leq h.
\end{align*}
Thus, every integer $n$ in the sumset $hA$ has a unique representation in the form 
$n = x_0a + y_0b$ with $x_0 \in [0,b-1]$ and $0 \leq y_0 \leq h - x_0$.
Conversely, if $n = x_0a+y_0b$ with $x_0 \in [0,b-1]$ and $y_0 \in [0,h-x_0]$, 
then $n$ is an element of the sumset $hA$.  Therefore,
\begin{align*}
|hA| & = \sum_{x_0=0}^{b-1} (h-x_0+1) = \frac{ b(2h+3-b)}{2} \\
& =  \binom{h+2}{2} - \binom{h+2-b}{2}. 
\end{align*}
Note that $\binom{h+2-b}{2} = 0$ for all $b \geq h+1$ 
and so it suffices to consider only sets $A = \{0,a,b\}$ with $0 < a < b \leq h+1$.
In particular, for all $b \in [2,h+1]$, the set $A = \{0,1,b\}$ has sumset 
\[
\binom{h+2}{2} - \binom{h+2-b}{2} =  \binom{h+2}{2} - \binom{t}{2}
\] 
for $t = h+2-b \in [1,h]$.
This completes the proof. 
\end{proof}

\bprob 
The structure of sumsets of sets of four integers is still not understood.  
Compute the set of sumset sizes $\mcr_{\Z}(h,4)$.
What is the cardinality of this set? 
\eprob 

\bprob
A general  problem is to compute $\mcr_{\Z}(h,k)$ for fixed $k$ 
as $h$ increases.
\eprob


\section{Sequences of sumset sizes}

Let $A$ be a   finite set of integers with $|A| = k \geq 2$.  
By translation, we can assume that 
$\min(A) = 0$ and so
\[
A \subseteq 2A \subseteq \cdots \subseteq hA \subseteq (h+1)A \subseteq \cdots. 
\]
If  $\max(A) = a > 0$, then  
\[
\max(hA) = ha < (h+1)a = \max ((h+1)A) 
\]
and so $hA$ is a proper subset of $(h+1)A$.  
Defining 
\[
k_h = |hA|
\] 
for all positive integers $h$, we obtain the strictly increasing sequence of positive integers 
\[
k = |A| = k_1 < k_2 < \cdots < k_h < \cdots.   
\]
Let 
\[
\kappa_{\infty}(A) = (k_1,k_2,\ldots, k_h, \ldots).  
\]
and, for every positive integer $\ell$, 
\[
\kappa_{\ell}(A) = (k_1,k_2,\ldots, k_{\ell}) 
\] 
By a fundamental theorem of additive number theory (Nathanson~\cite{nath1972-7}), 
the sequence $\kappa_{\infty}(A)$ is eventually an arithmetic progression 
with difference $\max(A)$. 

Consider the finite set 
\[
\mcr_{\Z}^{\sharp}(\ell,k) 
= \left\{ \kappa_{\ell}(A): A \subseteq \Z  \text{ and } |A| = k \right\}   
\]
and the infinite set 
\[
\mcr_{\Z}^{\sharp}(\infty,k) = \left\{ \kappa_{\infty}(A): A \subseteq \Z  \text{ and } |A| = k \right\}.     
\] 
For example, if $k=1$, then $k_h = |hA|=1$ for all $h$ and  
\[
\mcr_{\Z}^{\sharp}(\ell,1) = \left\{ (1,1,\ldots,1)  \right\} 
\]
for all $\ell$.  If $k=2$, then  $k_h = |hA|= h+1$ for all $h$  and   
\[
\mcr_{\Z}^{\sharp}(\ell,2) = \left\{ (2,3,\ldots, \ell+1)  \right\} 
\]
for all $\ell$. 
If $ k = 3$, then $k_2 \in \{5,6\}$, and $k_3 \in \{7,9,10\}$. 
We have 
\[
\kappa_3(0,1,2) = (3,5,7)
\]
\[
\kappa_3(0,1,3) = (3,6,9)
\]
\[
\kappa_3(0,1,4) = (3,6,10).
\]
Because $k_2=5$ if and only if $k_3=7$
 if and only if $A$ is an arithmetic progression, 
we obtain  
\[
\mcr_{\Z}^{\sharp}(3,3) = \left\{ (3,5,7), (3,6, 9 ), (3,6,10)  \right\}.
\]

\bprob
Compute $\mcr_{\Z}^{\sharp}(\ell,3)$ for $\ell \geq 4$.   
\eprob

\bprob
A general problem is to compute $\mcr_{\Z}^{\sharp}(\ell,k)$ for fixed $k$ 
as $\ell$ increases.
\eprob

\bprob
Thinking of the sequences $\kappa_{\ell}(A)$ and $\kappa_{\infty}(A)$ 
as trajectories of the sumset size function $k_h=|hA|$ 
through the space of positive integers suggests new problems.    
For example, given integers $k_1$ and $k_3$, determine all integers $k_2$ such that 
there exists a set $A$ with $|A|=k_1$, $|2A| = k_2$, and $|3A| = k_3$.  Equivalently, 
given  $k_1$ and $k_3$,
compute the set of ``intermediate states''
\[
\{k_2: (k_1,k_2,k_3) \in \mcr_{\Z}^{\sharp}(3,k_1) \}.
\] 
\eprob 

\bprob
Let $\ell \geq 4$.  Given integers $k_1$ and $k_{\ell}$, determine the  
sumset size trajectories from $k_1$ to $k_{\ell}$, 
that is, compute the set 
\[
\{ (k_2,k_3,\ldots, k_{{\ell}-1}): (k_1,k_2,k_3,\ldots, k_{\ell -1}, k_{\ell}) \in \mcr_{\Z}^{\sharp}({\ell},k_1) \}.
\] 
\eprob

\bprob
Let ${\ell} \geq 3$.  Given integers $k_1 < k_2 < \cdots <k_{\ell}-1$, 
compute the set of ``final states''
\[
\{k_{\ell}: (k_1,k_2,\ldots, k_{\ell -1},k_{\ell}) \in \mcr_{\Z}^{\sharp}(h,k_1) \}.
\] 
\eprob

\bprob
Here is a more general version of the previous problems.  
Let $h_1< \cdots < h_{p} \leq \ell$ and $k_1< \cdots < k_{p}$ 
be increasing sequences of positive integers.    
Consider the set 
\[
\mcr_{\Z}^{\flat} \left( (h_i)_{i = 1}^{p}, (k_i)_{i = 1}^{p} \right) 
= \left\{ \kappa_{\ell}(A): A \subseteq \Z \text{ and } |h_iA| = k_i \text{ for all } i \in [1, p] \right\}. 
\]
Find conditions on the sequences $(h_i)_{i=1}^{p}$  and $(k_i)_{i=1}^{p}$  
such that the set $\mcr_{\Z}^{\flat} \left( (h_i)_{i = 1}^{p}, (k_i)_{i = 1}^{p} \right)$ is nonempty.
Compute this set.
\eprob

\bprob 
There are unsolved  problems about the relation between sumsets and restricted
sumsets of a finite set of integers. For example, 
let $h$, $k$, $t$, and $\widehat{t}$ be positive integers such that 
\[
hk-h+1 \leq t \leq \binom{k+h-1}{h}
\]
and 
\[
hk-h^2+1 \leq \widehat{t} \leq \binom{k}{h}
\]
Decide if there exists a subset $A$ of $G$ such that 
\[
|A| = k, \qquad |hA| = t, \qquad  \left|\widehat{hA}\right| = \widehat{t}.  
\]
Compute the set  
\[
\left\{
\left( |hA|, \left|\widehat{hA}\right| \right): A \subseteq G \text{ and } |A| = k 
\right\}.
\]
\eprob


\appendix 

\section{On $|hA| = |A|$}         \label{sizes:appendix:hA=A}
\bl
Let $G$ be an additive abelian group and let $A$ be a nonempty finite subset of $G$.  
For all $h \geq 2$,  
\[
|hA| = |A|
\]
if and only if $A$ is a coset of a subgroup of order $|A|$.
\el

\begin{proof}
For all $x \in G$, the translation of the set $A$ by $x$ is the set
\[
x+A = \{x+a:a\in A\}.
\]
For all $h \geq 2$ we have 
\[
h(x+A) = hx+hA 
\]
and so 
\[
|h(x+A) |  = |hA| 
\]
By translation, we can assume that $0 \in A$ and so 
\[
A \subseteq 2A \subseteq 3A \subseteq \cdots \subseteq hA. 
\]
It follows that 
\[
|A| \leq |2A| \leq |3A| \leq \cdots \leq |hA| .
\]
If $|hA| = |A|$, then 
\[
A = 2A =3A = \cdots = hA.
\]
In particular, $A = 2A$ and so $A$ is closed under addition in the group $G$.
For all $x \in A$, we have 
\[
x+A \subseteq A+A = A.
\]
Because $|x+A| = |A|$, we have $x+A=A$ and so there exists $y \in A$ 
such that $x+y=0$. 
Thus, $A$ is a subgroup of $G$.  
This completes the proof. 
\end{proof}

\section*{Acknowledgements}

I thank Kevin O'Bryant for bringing the Erd\H os-Szemer\' edi  statement to my attention
  and sparking an interest in the sizes of sumsets.  

\end{document}